\documentclass[12pt]{article}
\usepackage{amsfonts,amscd,amssymb}
\input diagrams

\newtheorem{lemma}{Lemma}[section]
\newtheorem{theorem}[lemma]{Theorem}
\newtheorem{proposition}[lemma]{Proposition}
\newtheorem{corollary}[lemma]{Corollary}

\newtheorem{definition}[lemma]{Definition}

\newtheorem{example}[lemma]{Example}

\newenvironment{proof}{{\it Proof.}}{\hfill $ \square $ \vskip 4mm}

\newcommand{\eqref}[1]{(\ref{eq:#1})}

\def\sw#1{{\sb{(#1)}}}
\def\sco#1{{\sp{(\bar #1)}}} 
\def\su#1{{\sp{(#1)}}} 

\def\<{{\langle}}
\def\>{{\rangle}}

\def\eps{\epsilon}

\def\note#1{{}}
\def\can{{\rm can}}
\def\note#1{}

\def\cocan{{\rm cocan}}

\def\B{{\rm Bar}}
\def\C{{\rm Cob}}

\begin{document}
\title{The cohomology structure of an algebra 
entwined with a coalgebra}
\author{Tomasz Brzezi\'nski
\thanks{Lloyd's of London Tercentenary Fellow at the University of York.}\\
Department of Mathematics, University of York\\
 Heslington
York YO10 5DD, England\thanks{Address for correspondence. E-mail:
tb10@york.ac.uk} \\ \&\\ 
Department of Theoretical Physics, University of \L\'od\'z\\ Pomorska
149/153, 90-236
\L\'od\'z, Poland}
\date{}
\maketitle

\begin{abstract}
Two cochain complexes are constructed for an algebra $A$ and 
 a coalgebra $C$ entwined with each other via the map $\psi:C\otimes
A\to A\otimes C$. One complex is associated to an $A$-bimodule, 
the other to a
$C$-bicomodule. In the former case the resulting complex can be
considered as a $\psi$-twisted Hochschild complex of $A$, while for the
latter one obtains a $\psi$-twist of the Cartier complex of $C$. 
The notion of a {\em weak comp algebra} is introduced
by weakening the axioms of the Gerstenhaber comp algebra. It is shown
that such a  weak comp algebra is a cochain complex with two cup
products that descend to the cohomology. It is also shown that the
complexes associated to an entwining structure and $A$ or $C$ are
examples of a weak comp algebra.
Finally both complexes are combined in a double complex whose role in
the deformation theory of entwining structures is outlined.
\end{abstract}
\baselineskip=19pt

\section{Introduction}
An entwining structure \cite{BrzMaj:coa} comprises of an algebra, 
a coalgebra and a map
that entwines one with the other and satisfies some simple axioms.
In many respects an entwining structure
resembles a bialgebra or a comodule algebra of a bialgebra. Indeed, to
any comodule algebra of a bialgebra, and hence to a bialgebra itself, 
 there is associated a generic
entwining structure, canonical in a certain sense. The aim of this
paper is to reveal that entwining structures admit a rich cohomology 
theory,
which, depending on the choice of ingredients in the entwining
structure, can be viewed as the Hochschild cohomology of an algebra
\cite{Hoc:coh} or 
the Cartier cohomology of a coalgebra \cite{Car:coh}, and is reminiscent
of the Gerstenhaber-Schack theory for bialgebras \cite{GerSch:bia}.

Recall from \cite{BrzMaj:coa} that an {\em entwining structure} over a field
$k$ consists of an algebra $A$, a coalgebra $C$ and a map $\psi:C\otimes
A\to A\otimes C$ such that the following {\em bow-tie diagram} commutes.
$$
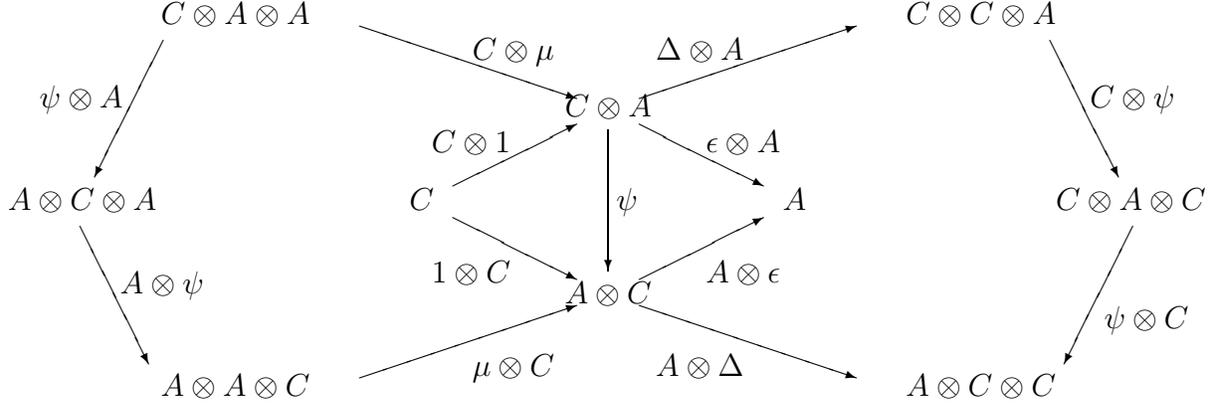
\begin{diagram}
& C\otimes A\otimes A &  & & & & & & &
C\otimes C\otimes A &  \\
\ldTo(1,2)^{\psi\otimes A} &  & & \rdTo(3,1)^{C\otimes\mu} & &  C\otimes
A & & \ruTo(3,1)^{\Delta\otimes A} &&& \rdTo(1,2)^{C\otimes\psi} \\
A\otimes C\otimes A &&&  C& \ruTo(2,1)^{C\otimes 1} & \dTo_\psi &
\rdTo(2,1)^{\eps\otimes A}& A &&& C\otimes A\otimes C \\
\hspace{2.1cm} \rdTo(1,2)^{A\otimes\psi} &&&& \rdTo(2,1)_{1\otimes C} 
&A\otimes C & 
\ruTo(2,1)_{A\otimes \eps} &&&& \hspace{-2.1cm} \ldTo(1,2)_{\psi\otimes C}\\
&A\otimes A\otimes C &&&\ruTo(3,1)_{\mu\otimes C}
&&\rdTo(3,1)_{A\otimes\Delta}  &&& A\otimes C\otimes C &
\end{diagram}
$$

Here and below we use the following notation. The product in $A$ is
denoted by $\mu$, while the unit (both as an element of $A$ and the
map $k\to A$) is denoted by $1$. For a coalgebra $C$, $\Delta$ is
the coproduct, while $\eps$ is the counit. We use the Sweedler
notation for action of $\Delta$ on elements of $C$, $\Delta(c) =
c\sw 1\otimes c\sw 2$ (summation understood). Finally for any vector
space $V$, the identity map $V\to V$ is denoted by $V$; we also
implicitly identify $k\otimes V$ and $V\otimes k$ with $V$, use
$\otimes$ for $\otimes_k$ and write $V^n$ for $V^{\otimes n}$.

An entwining
structure is denoted by $(A,C)_\psi$. To describe the action of $\psi$ 
we  use the following {\it $\alpha$-notation}: 
$\psi(c\otimes a) = a_\alpha \otimes c^\alpha$ (summation
over a Greek index understood), for all $a\in
A$, $c\in C$, which proves very useful in concrete computations
involving $\psi$. Reader is advised to check that the bow-tie diagram is
equivalent to the following four explicit relations:
\begin{eqnarray*}
&&\hspace{-.7cm} \mbox{left pentagon:}\;\; (aa')_\alpha \otimes c^\alpha = 
a_\alpha a'_\beta \otimes
c^{\alpha\beta} , \quad \mbox{left triangle:}\;\;
1_{\alpha}\otimes c^\alpha = 1\otimes c,\\
&&\hspace{-.7cm}\mbox{right pentagon:}\; 
a_\alpha\otimes c^\alpha\sw 1\otimes c^\alpha\sw 2 =
a_{\beta\alpha}\otimes 
c\sw 1^\alpha\otimes c\sw 2^\beta , \quad 
\mbox{right triangle:}\; a_\alpha\eps(c^\alpha) = 
a\eps(c),
\end{eqnarray*}
for all $a,a'\in A$, $c\in C$.  

One may (or perhaps even should) think of an entwining map $\psi$ as a
twist in the convolution algebra ${\rm Hom}(C,A)$. Namely, given an
entwining structure, one can define  the map $*_\psi: 
{\rm Hom}(C,A)\otimes{\rm Hom}(C,A)\to{\rm Hom}(C,A)$ via $(f*_\psi
g)(c) = f(c\sw 2)_\alpha g(c\sw 1^\alpha)$, for all $f,g\in {\rm
Hom}(C,A)$ and $c\in C$. One can easily check that $({\rm
Hom}(C,A),*_\psi)$ is an associative algebra with unit $1\circ\eps$.

There are many examples of entwining structures. As a generic example
one can refer  to the following situation. Suppose $C$ is a bialgebra
and  $A$ is a right $C$-comodule algebra with the coaction $\rho^A$. 
Then $\psi:c\otimes a \mapsto a\sw 0 \otimes ca\sw 1$, where $\rho^A(a) = a\sw
0\otimes a\sw 1$, entwines $C$ with $A$. A special case of this
situation is when $A$ is an algebra and a coalgebra at the same time.
Then $A$ is a bialgebra if and only if $\psi:A\otimes A \to A\otimes A$,
$\psi:a\otimes a'\mapsto a'\sw 1\otimes a a'\sw 2$ entwines $A$ with
itself. Furthermore, any algebra and a
coalgebra can be provided with an entwining structure with $\psi$ being
the usual flip of tensor factors (for obvious reasons this can be called
a {\em trivial} entwining structure). Interesting examples come from
the generalisation of Hopf-Galois theory, motivated by the geometry of
quantum (group) homogeneous spaces.

\begin{example}\rm\cite{BrzHaj:coa}
Let $C$ be a coalgebra, $A$ an algebra and a right $C$-comodule with the
coaction $\rho^A:A\to A\otimes C$. Let
$B:= \{b\in A\; | \; \rho^A(ba) =
b\rho^A(a)\}$ and assume that the canonical left $A$-module, right
$C$-comodule map
$
\can:\ A\otimes _BA\to A\otimes C$, $a\otimes a'\mapsto a\rho^A(a')$,
is bijective. Let $\psi:C\otimes A\to A\otimes C$ be a $k$-linear map
given by
$\psi(c\otimes a) = \can(\can^{-1}(1\otimes c)a)$.
Then $(A,C)_\psi $ is an entwining structure. The
extension $B\hookrightarrow A$ is called a
{\em coalgebra-Galois extension} (or a {\em $C$-Galois extension})  
and is
denoted by $A(B)^C$. This is a generalisation of the notion of 
a Hopf-Galois extension
introduced in \cite{KreTak:hop}. $(A,C)_\psi$ is called a
{\em canonical entwining structure} associated to $A(B)^C$. 
\label{can.ex}
\end{example}

If $C$ is a Hopf algebra and $A(B)^C$ is a Hopf-Galois extension, then
the canonical entwining structure is the generic one described above.
Also, any Hopf algebra is a Hopf-Galois extension of $k$, and the
canonical entwining structure in this case is the one described above
for a bialgebra. Dually we have

\begin{example}\rm\cite{BrzHaj:coa}
Let $A$ be an algebra, $C$ a coalgebra and a right $A$-module with the
action $\rho_C: A\otimes C\to C$. Let
$B:= C/I$, where $I$ is a coideal in $C$,
$$
I:={\rm span}\{(c\cdot a)\sw 1\xi((c\cdot a)\sw 2)-c\sw 1\xi(c\sw
2\cdot a)\; |\; a\in A,\; c\in C,\; \xi\in C^*\},
$$
 and assume that the canonical left $C$-comodule, right
$A$-module map
$
\cocan:\ C\otimes A\to C\square_B C$, $c\otimes a \mapsto c\sw 1\otimes c\sw
2\cdot a$,
is bijective. Let $\psi:C\otimes A\to A\otimes C$ be a $k$-linear map
given by
$$
\psi = (\eps\otimes A\otimes C)\circ (\cocan^{-1}\otimes C)\circ
(C\otimes \Delta)\circ \cocan.
$$
Then $(A,C)_\psi $ is an entwining structure. The
coextension $C\twoheadrightarrow B$ is called an
{\em algebra-Galois coextension} (or an {\em $A$-Galois coextension})
and is
denoted by $C(B)_A$. $(A,C)_\psi$ is called a
{\em canonical entwining structure} associated to $C(B)_A$. 
\label{cocan.ex}
\end{example}

In dealing with cohomology we will need $A$-bimodule
($C$-bicomodule resp.) structures on $A\otimes C^n$ ($C\otimes A^n$
resp.). These are defined as follows. 
Given an entwining structure $(A,C)_\psi$ define two infinite 
families of maps
$$
\psi^n  = (A^{n-1}\otimes \psi)\circ (A^{n-2}\otimes\psi\otimes
A)\circ\cdots \circ (\psi\otimes A^{n-1}) : 
C\otimes A^n\to A^n\otimes C, \quad  n\geq 1, 
$$
$$
\psi_n = (\psi\otimes
C^{n-1})\circ(C\otimes\psi\otimes C^{n-2})\circ\cdots\circ
(C^{n-1}\otimes \psi) : C^n\otimes A \to A\otimes C^n, \quad
n\geq 1.
$$
The axioms of an entwining structure imply that for all $n>0$, 
$A\otimes C^n$ is an
$A$-bimodule with the left action $\rho_n^L=\mu\otimes
C^n$ and the right action $\rho_n^R = (\mu\otimes C^n)\circ
(A\otimes \psi_{n})$. Furthermore,   $C\otimes A^n$ is a
$C$-bicomodule with the left coaction $\rho^n_L=\Delta
\otimes A^n$ and the right coaction $\rho^n_R = 
(C\otimes \psi_{n})\circ (\Delta
\otimes A^n)$. We will always consider $A\otimes C^n$ ($C\otimes
A^n$ resp.) as bimodules (bicomodules resp.) with the above structures.
Also, for any vector space $V$ the space $A\otimes V\otimes A$
($C\otimes V\otimes C$ resp.) will be considered as an $A$-bimodule
($C$-bicomodule resp.) with the obvoius structure maps $\mu\otimes V\otimes A$ and
$A\otimes V\otimes \mu$ ($\Delta\otimes V\otimes C$, $C\otimes V\otimes
\Delta$ resp.).

Yet another consequence of  the axioms of an entwining structure 
is the following 
\begin{lemma}
For all $n\geq 1$, $0\leq j\leq n-1$ the following two diagrams
\begin{diagram}
C\otimes A^{n+1} & \rTo^{C\otimes A^j\otimes \mu\otimes A^{n-j-1}} &
C\otimes A^{n}  \\
\dTo^{\rho^{n+1}_R} & & \dTo_{\rho^n_R}  \\
C\otimes A^{n+1}\otimes & \rTo^{C\otimes A^j\otimes\mu\otimes
A^{n-j-1} \otimes C}&  C\otimes A^n\otimes C 
\end{diagram}
\begin{diagram}
 A\otimes C^n\otimes A & 
\rTo^{A\otimes C^j\otimes \Delta\otimes C^{n-j-1}
\otimes A} & A\otimes C^{n+1}\otimes A \\
\dTo^{\rho_n^R} & & \dTo_{\rho_{n+1}^R} \\
A\otimes C^{n} & \rTo^{A\otimes C^j\otimes\Delta\otimes
C^{n-j-1}} &  A\otimes C^{n+1}  
\end{diagram}
commute.
\label{lemma.system}
\end{lemma}

The paper is organised as follows. In the next section we construct the
cochain complex $C_\psi(A,M)$ associated to an entwining structure $(A,C)_\psi$ and an
$A$-bimodule
$M$. We study its relation to the Hochschild complex of $A$ as well as
analyse its structure in the case of the canonical entwining structure
associated to a $C$-Galois extension $A(B)^C$. In particular we show
that if $B=k$ this complex provides a resolution of $M$. In Section~3 we
dualise the construction of Section~2 and describe a complex
$A_\psi(C,V)$ associated to
$(A,C)_\psi$ and a $C$-bicomodule $V$. Section~4 is devoted to
studies of the complex $C_\psi(A) = C_\psi(A,A)$. We
define two cup products in $C_\psi(A)$ which descend to the 
cohomology. In particular we show
that in the cohomology one product is a graded twist of the other. All
this is done with the help of the notion of a {\em  weak comp
algebra}, which generalises the notion of a right comp algebra
\cite{GerSch:alg} or a pre-Lie system \cite{Ger:coh} (see recent review
\cite{Klu:inv} of comp algebras). In
Section~5 we
define an {\em equivariant} complex as a subcomplex of
$C_\psi(A)$ on which both cup products coincide, so that the
corresponding algebra in such an {\em equivariant} 
cohomolgy is graded commutative. This extends the classic result of
Gerstenhaber \cite{Ger:coh}. Finally in Section~6 we define a double
complex and outline its role in the deformation theory of entwining
structures. 

\section{Module valued cohomology of  an entwining structure}
In this section we associate a cochain complex 
to an entwining structure $(A,C)_\psi$ and an $A$-bimodule $M$.

Recall that the bar resolution of an  algebra $A$ is a chain
complex $\B(A)= (\B_\bullet(A),\delta)$, where
$$
\B_n(A) = A^{n+2},\quad \delta_n =
\sum_{k=0}^{n} (-1)^k A^{k}\otimes\mu\otimes A^{n-k} :
A^{n+2}\to A^{n+1}.
$$
Since $\delta$ is an $A$-bimodule map, one can define
a chain complex $\B^\psi(A)= (\B_\bullet^\psi(A),\delta)$ via
$
\B^\psi(A)= (A\otimes C)\otimes_A\B(A).
$
Here $A\otimes C$ is viewed as an $A$-bimodule as explained in the
introduction.  Explicitly, $\B_n^\psi(A) = A\otimes C\otimes A^{n+1}$ and 
$$
\delta_n =
(\mu\otimes C\otimes A^n)\circ(A\otimes\psi\otimes A^n)
+\sum_{k=1}^n(-1)^k(A\otimes C\otimes A^{k-1}\otimes\mu\otimes
A^{n-k}).
$$
Even more explicitly, using the $\alpha$-notation, one can write
\begin{eqnarray*}
\delta_n(a^0\otimes c\otimes a^1\otimes\cdots\otimes a^{n+1}) & = &
a^0a^1_\alpha\otimes c^\alpha\otimes a^2\otimes\cdots\otimes a^{n+1}\\
&&+ \sum_{i=1}^{n}(-1)^ia^0\otimes c\otimes a^1\otimes\cdots\otimes 
a^ia^{i+1}\otimes\cdots\otimes a^{n+1}.
\end{eqnarray*}

\begin{lemma}
$\B^\psi(A)$ 
is a resolution of both $A\otimes C$ and $A$. Furthermore, $\delta$ is an
$A$-bimodule map.
\end{lemma}
\begin{proof}
$\B^\psi(A)$ is an acyclic complex, because there is a 
contracting
homotopy $h_n: \B_n^\psi(A)\to  \B_{n+1}^\psi(A)$, given as $h_n = 
(-1)^n A\otimes
C\otimes A^{n+1}\otimes 1$. The left pentagon in the bow-tie diagram
implies that the map  $(\mu
\otimes C)\circ (A\otimes\psi):A\otimes C\otimes A\to A\otimes C$ is an
augmentation. Similarly, the right triangle implies that
 $\mu\circ(A\otimes\eps\otimes A): A\otimes C\otimes A\to A$ is an
augmentation as well.
 It is obvious that all the $\delta_n$ are $A$-bimodule
maps.
\end{proof}

Next we use the resolution $\B^\psi(A)$ to construct 
one of the main 
cochain complexes studied in this paper. 
Let $M$ be an $A$-bimodule. Define the cochain complex $C_\psi (A,M) = 
(C_\psi(A,M)^\bullet, d)$ by $ C_\psi (A,M) ={}_A{\rm
Hom}_A(\B^\psi(A), M)$, where ${}_A{\rm Hom}_A$ denotes the Hom-bifunctor
from the category of 
$A$-bimodules to the category of $k$-vector spaces.
Using the natural identification
$
_A{\rm
Hom}_A(A\otimes C\otimes A^{n+1}, M) = {\rm Hom}(C\otimes A^n,M),
$
one explicitly obtains
$$
C^n_\psi(A,M) = {\rm Hom}(C\otimes A^n,M), \quad d^n : {\rm Hom}(C\otimes
A^n,M) \to  {\rm Hom}(C\otimes A^{n+1},M),
$$
$$
 d^n f =
{}_M\rho\circ(A\otimes f)\circ(\psi\otimes A^n) +\sum_{k=1}^n(-1)^k f\circ
(C\otimes A^{k-1}\otimes\mu\otimes A^{n-k})
+(-1)^{n+1}\rho_M\circ(f\otimes A),
$$
where $\rho_M$ (${}_M\rho$ resp.) denotes the right (left) action of $A$
on $M$. Even more explicitly, using the $\alpha$-notation, one has
\begin{eqnarray*}
d^{n}f(c,a^1,\ldots ,a^{n+1}) & = & a^1_\alpha\cdot 
f(c^\alpha,a^2,\ldots ,a^{n+1})
+ \sum_{i=1}^n (-1)^i 
f(c,a^1,\ldots ,a^ia^{i+1},\ldots ,a^{n+1})\\
&&+ (-1)^{n+1}f(c,a^1,\ldots ,a^{n})\cdot a^{n+1}.
\end{eqnarray*}
To save the space we write $f(\cdot, \cdot, \ldots, \cdot)$ for
$f(\cdot\otimes\cdot\otimes\cdots\otimes\cdot)$, etc. 
There is a close relationship between the complex 
$C_\psi(A, M)$ and the Hochschild complex of $A$.
Firstly, for any $A$-bimodule $M$ consider the
$A$-bimodule ${\rm Hom}(C,M)$ with the structure maps:
$$
(f\cdot a)(c) = f(c)\cdot a, \qquad (a\cdot f)(c) = a_\alpha\cdot 
f(c^\alpha).
$$
It is an easy exercise in the $\alpha$-notation to verify that ${\rm
Hom}(C,M)$ is an $A$-bimodule indeed. 
Identify ${\rm Hom}(A^n,{\rm Hom}(C,M))$ with 
${\rm Hom}(C\otimes A^n,M)$ by the natural isomorphism:
$$
\theta(f)(c,a^1,\ldots ,a^n) = f(c)(a^1,\ldots ,
a^n), \quad \theta^{-1}(g)(a^1,\ldots , a^n)(c) = g(c, a^1,\ldots ,a^n).
$$
Then $C_\psi(A, M)$ is  the Hochschild complex of
$A$ with values in ${\rm Hom}(C,M)$.

Secondly, the Hochschild complex over $A$ with values in
$M$ is included in the complex $C_\psi(A,M)$. More
precisely one has
\begin{lemma}
Let $C(A,M)$ be the Hochschild complex over
$A$ with values in $M$. 
Then the map $j:C(A, M) \to C_\psi(A, M)$ given
by 
$$
j^n :{\rm Hom}(A^n,M)\to {\rm Hom}(C\otimes A^n, M), \qquad j^n:
f\mapsto \eps\otimes f,
$$
is a monomorphism of differential complexes.
\label{lemma.hoch}
\end{lemma}
\begin{proof}
Clearly, $j$ is injective. The fact that $j$ is the map between cochain
complexes follows from the right triangle in the bow-tie diagram.
\end{proof}

The cohomology of the complex $C_\psi(A, M)$ is
denoted by $H_\psi(A, M)$ and is called an {\em entwined
cohomology of $A$ with values in $M$}. 

\begin{proposition}
For an entwining structure $(A,C)_\psi$, $A\otimes C$ is a projective $A$-bimodule if and only if
$H^1_\psi(A,M)=0$ for all $A$-bimodules $M$.
\label{prop.proj}
\end{proposition}

Proposition~\ref{prop.proj} will follow from the following two lemmas.
\begin{lemma}
For an entwining structure $(A,C)_\psi$ the following statements are
equivalent: 

(1) $A\otimes C$ is a projective $A$-bimodule.

(2) The sequence of $A$-bimodule maps
\begin{diagram}
 0 & \rTo & \ker{\rho_1^R} &\rTo & A\otimes C\otimes A & 
\rTo^{\rho_1^R} & A\otimes C &\rTo & 0
\end{diagram}
is split exact.

(3) There exists a $0$-cocycle $\chi\in C^0_\psi(A,A\otimes A)$ such
that $\mu\circ\chi = 1\circ\eps$.
\label{lemma.proj1}
\end{lemma}
\begin{proof}
The equivalence of the first two assertions is clear since $k$ is the
field so that $A\otimes
C\otimes A$ is a free $A$-bimodule. Suppose that (2) holds, i.e.,
there is an $A$-bimodule map $\nu:A\otimes C\to A\otimes C\otimes A$
such that $\rho^R_1\circ\nu = A\otimes C$.  Define $\bar{\chi} =
\nu\circ (1\otimes C)$ and $\chi= (A\otimes\eps\otimes
A)\circ\bar{\chi}$. Then for all $a\in A$, $c\in C$ one has:
$$
a_\alpha\bar{\chi}(c^\alpha) = a_\alpha\nu(1\otimes c^\alpha) =
\nu((1\otimes c)\cdot a) = \nu(1\otimes c)a = \bar{\chi}(c)a.
$$
Applying $A\otimes\eps\otimes A$ to both sides of this equality one
immediately obtains that $\chi$ is a $0$-cocycle. The normalisation of
$\chi$ follows from the equality $\rho^R_1\circ\nu = A\otimes
C$ applied to $1\otimes c$ and the right triangle in the bow-tie
diagram.

Now suppose that (3) holds. Denote $\chi(c) = c\sco 1\otimes c\sco 2$
(summation understood), and define $\nu:A\otimes C\to A\otimes C\otimes
A$, $a\otimes c\mapsto ac\sw 2\sco 1_\alpha\otimes c\sw 1^\alpha\otimes
c\sw 2\sco 2$. Then for all $a\in A$, $c\in C$,
$$
\rho^R_1\circ\nu(a\otimes c) = ac\sw 2\sco 1_\alpha c\sw 2\sco 2_\beta 
\otimes c\sw 1^{\alpha\beta} = a(c\sw 2\sco 1 c\sw 2\sco 2)_\alpha 
\otimes c\sw 1^{\alpha} = a\otimes c,
$$
where we used the left pentagon and the left triangle in the bow-tie
diagram together with the normalisation of $\chi$. Therefore $\nu$
splits $\rho_1^R$. Clearly, $\nu$ is a left $A$-module map. Furthermore
for all $a,a'\in A$, $c\in C$ we have:
\begin{eqnarray*}
\nu((a\otimes c)\cdot a') &=& aa'_\beta c^\beta\sw 2\sco 1_\alpha
\otimes c^\beta\sw 1^\alpha\otimes
c^\beta\sw 2\sco 2\\
&=& aa'_{\beta\gamma} {c\sw 2^{\beta\sco 1}}_\alpha
\otimes c\sw 1^{\gamma\alpha}\otimes
c\sw 2^{\beta\sco 2}\\
&=& a(a'_\beta c\sw 2^{\beta\sco 1})_\alpha
\otimes c\sw 1^\alpha\otimes
c\sw 2^{\beta\sco 2}\\
&=&ac\sw 2\sco 1_\alpha\otimes c\sw 1^\alpha\otimes
c\sw 2\sco 2a' = \nu(c\otimes a)a',
\end{eqnarray*}
where we used the right pentagon to derive the second equality, the left
pentagon to derive the third one and finally the fact that $\chi$ is a
$0$-cocycle to obtain the fourth equality. This proves that $\nu$ is an
$A$-bimodule splitting as required.
\end{proof}

\begin{lemma}
For an entwining structure $(A,C)_\psi$ and an $A$-bimodule $M$ let
$B^n_\psi (A,M)$ denote the space of $n$-coboundaries and 
$Z_\psi^n(A,M)$ denote the space of $n$-cocycles in $C_\psi^n(A,M)$. Let
$D_\psi :C\otimes A\to A\otimes C\otimes A$ be a linear map given by
$D_\psi : c\otimes a \mapsto 1\otimes c\otimes a - a_\alpha\otimes
c^\alpha\otimes 1$. Then:

(1) The map $\theta :{}_A{\rm Hom}_A(\ker\rho_1^R, M) \to
Z_\psi^1(A,M)$, $f\mapsto f\circ D_\psi$ is a bijection.

(2) $\theta^{-1}(B^n_\psi (A,M)) = \{ f\mid_{\ker\rho^R_1} \; :\; f\in 
{}_A{\rm Hom}_A(A\otimes C\otimes A, M)\}$.
\label{lemma.proj2}
\end{lemma}
\begin{proof}
Throughout the proof of this lemma, 
$x=\sum_i a^i\otimes c_i\otimes\tilde{a}^i$ is an arbitrary element
of $\ker\rho^R_1$. Notice that $\sum_i a^i\tilde{a}^i_\alpha\otimes
c_i^\alpha =0$.

(1) One easily finds that for all $f\in {}_A{\rm Hom}_A(\ker\rho_1^R, M)$, 
$d\theta(f) =0$ so that the map $\theta$ is well-defined. Consider
the map $\bar{\theta} :Z_\psi^1(A,M)\to {}_A{\rm Hom}_A(\ker\rho_1^R, M)$,
given by $\bar{\theta}(\chi)(x)
= \sum_i a^i\cdot\chi(c_i,\tilde{a}^i)$. Clearly, $\bar{\theta}(\chi)$ is a left
$A$-module map. Since $\chi$ is a 1-cocycle, 
we have for all $a\in A$,
$$
\bar{\theta}(\chi)(xa)= \sum_i a^i\cdot\chi(c_i,\tilde{a}^ia) =
\sum_i a^i\cdot\chi(c_i,\tilde{a}^i)a + \sum_i a^i\tilde{a}^i_\alpha
\cdot\chi(c_i^\alpha,a) = \bar{\theta}(\chi)(x)\cdot a.
$$
Therefore $\bar{\theta}$ is well-defined. For any $\chi\in Z_\psi^1(A,M)$
one easily finds that for all $c\in C$, $\chi(c,1)=0$. Using this fact
one obtains
$$
\theta\circ\bar{\theta}(\chi)(c,a) = \bar{\theta}(\chi)(1\otimes c\otimes a
-a_\alpha\otimes c^\alpha\otimes 1) = \chi(c,a) - a_\alpha\chi(c^\alpha,1)
=\chi(c,a),
$$
as well as
$$
\bar{\theta}\circ\theta(f)(x) = \sum_i a^i\theta(f)(c_i,\tilde{a}^i)
= f(x-\sum_ia^i\tilde{a}^i_\alpha\otimes c^\alpha\otimes 1) = f(x).
$$
Therefore $\bar{\theta}$ is the inverse of $\theta$, $\bar{\theta} = 
\theta^{-1}$.

(2) Suppose $\chi = -df$ for some $f\in {\rm Hom}(C,M)$. Then
$\theta^{-1}(\chi)(x) = -\sum_ia^i\cdot df(c_i,\tilde{a}^i) = 
\sum_i a^i\cdot f(c_i)\cdot\tilde{a}^i - \sum_i a^i\tilde{a}^i_\alpha
\cdot f(c_i^\alpha) = \sum_i a^i\cdot f(c_i)\cdot \tilde{a}^i$. The result then follows
from the isomorphism ${\rm Hom}(C,M) \cong  {}_A{\rm Hom}_A(A\otimes A, 
{\rm Hom}(C,M)) \cong  {}_A{\rm Hom}_A(A\otimes C\otimes A, M)$ given by
$f\mapsto \ell_f$, where $\ell_f(a\otimes c\otimes a') = a\cdot f(c)\cdot
a'$.
\end{proof}

Now, using Lemma~\ref{lemma.proj1} and Lemma~\ref{lemma.proj2}, 
Proposition~\ref{prop.proj} can be proven by the same reasoning as
Proposition~11.5 in \cite{Pie:ass}. Namely, if $H^1_\psi(A,\ker\rho_1^R)=0$,
then, by Lemma~\ref{lemma.proj2}, $A$-bimodule endomorphisms of 
$\ker\rho_1^R$ equal $\{ f\mid_{\ker\rho^R_1} \; :\; f\in 
{}_A{\rm Hom}_A(A\otimes C\otimes A, \ker\rho^R_1)$. This means
that there exists $f\in 
{}_A{\rm Hom}_A(A\otimes C\otimes A, \ker\rho^R_1)\}$ such that 
$f\mid_{\ker\rho^R_1} = \ker\rho^R_1$, i.e., the sequence in 
Lemma~\ref{lemma.proj1}(2) splits. Thus $A\otimes C$
is a projective $A$-bimodule,  by Lemma~\ref{lemma.proj1}.
Conversely, if there is an extension $f\in 
{}_A{\rm Hom}_A(A\otimes C\otimes A, \ker\rho^R_1)$
of the identity mapping $\ker\rho_1^R$, then for
any $A$-bimodule $M$, every
$g\in 
{}_A{\rm Hom}_A(\ker\rho^R_1,M)$ has the form $h\mid_{\ker\rho^R_1}$,
where $h=g\circ f\in {}_A{\rm Hom}_A(A\otimes C\otimes A, M)$. By
Lemma~\ref{lemma.proj2} $B^1_\psi (A,M)=Z_\psi^1(A,M)$, for any 
$A$-bimodule $M$. This completes the proof of Proposition~\ref{prop.proj}.

As an example we compute the (canonical) entwined cohomology of a 
$C$-Galois extension.
\begin{proposition}
Let $(A,C)_\psi$ be the canonical entwining structure associated to a
$C$-Galois extension $A(B)^C$ in Example~\ref{can.ex}, and let $M$
be an $A$-bimodule. Then

(1)  $H^0_\psi(A,M) = M^B := \{m\in M\; | \; \forall b\in B, \; b\cdot m
= m\cdot b\}$.

(2) If $A$ is a $C$-Galois object, i.e., $B=k$, then $H^n_\psi(A,M) =
0$, for all $n>0$.
\label{prop.Galois}
\end{proposition}
\begin{proof}
(1) Denote the action of the {\em translation map}
 $\tau := \can^{-1}\circ (1\otimes C):
C\to A\otimes_B A$ by $\tau(c) = c\su 1\otimes c\su 2$ (summation
understood). Notice that for all $a\in A$, $c\in C$, 
$c\su 1c\su 2\sw 0\otimes c\su 2\sw 1 =
1\otimes c$, $a\sw 0a\sw 1\su 1\otimes a\sw 1\su 2 = 1\otimes a$ (cf.\
\cite[3.4~Remark~(2)(a)]{Sch:rep}), 
and $\psi(c\otimes a) = c\su 1(c\su 2a)\sw 0\otimes (c\su
2a)\sw 1$, where we use the Sweedler notation for the coaction,
$\rho^A(a) = a\sw 0\otimes a\sw 1$ (summation understood).
 Consider the map
$$
\theta: M^B\to H^0_\psi(A,M) , \qquad m\mapsto c\su 1\cdot m\cdot c\su
2.
$$
The map $\theta$ is well-defined because $m\in M^B$ and furthermore,
\begin{eqnarray*}
d\theta(m)(c\otimes a) &=& a_\alpha\cdot \theta(m)(c^\alpha) -
\theta(m)(c)\cdot a\\
&=& c\su 1(c\su 2a)\sw 0(c\su 2a)\sw 1\su 1\cdot m\cdot (c\su 2a)\sw 1\su 2 -
c\su 1\cdot m\cdot c\su 2a\\
&=& c\su 1\cdot m\cdot c\su 2a - c\su 1\cdot m\cdot c\su 2a =0,
\end{eqnarray*}
so that $\theta(m)$ is a zero-cocycle, hence belongs to the
zero-cohomology group.

Let $\rho^A(1) = 1\sw 0\otimes 1\sw 1$. We claim that the map
$$
\theta^{-1}: H^0_\psi(A,M)\to M^B, \qquad \theta^{-1}(f) = 1\sw 0\cdot
f(1\sw 1)
$$
is the inverse of $\theta$. Firstly we need to check that
$\theta^{-1}(f)$ is in the centraliser of $B$ in $M$. The key
observation needed here is that $A$ is an $(A,C)_\psi$-module, i.e., for
all $a,a'\in A$, $\rho^A(aa') = a\sw 0 a'_\alpha\otimes a\sw 1^\alpha$. 
In particular this implies that for
all $a\in A$, 
$\rho^A(a) = 1\sw 0 a_\alpha\otimes 1\sw 1^\alpha$, and for
all $b\in B$, $b1\sw 0\otimes 1\sw 1 = 1\sw 0 b_\alpha\otimes
1\sw 1^\alpha$ (cf.\ \cite{Brz:mod}). 
Since $f$ is a
zero-cocycle we have for any $b\in B$:
$$
0 = 1\sw 0\cdot df(1\sw 1,b) = 1\sw 0b_\alpha\cdot f(1\sw 1^\alpha) -
1\sw 0\cdot f(1\sw 1)\cdot b = b\cdot \theta^{-1}(f) - 
\theta^{-1}(f)\cdot b ,
$$
so that $\theta^{-1}(f)\in M^B$ as claimed. Furthermore, since
$\rho^A(a) = 1\sw 0a_\alpha\otimes 1\sw 1^\alpha$, we have for all
0-cocycles $f$, $a\sw 0\cdot f(a\sw 1) = 1\sw 0\cdot f(1\sw 1)\cdot a$.
In particular we have
$$
f(c) = c\su 1c\su 2\sw 0\cdot f(c\su 2\sw 1) = c\su 1 1\sw 0\cdot
f(1\sw 1)\cdot c\su 2.
$$
Therefore
$$
\theta\circ\theta^{-1}(f)(c) = c\su 1\cdot\theta^{-1}(f)\cdot c\su 2 =
 c\su 1 1\sw 0 \cdot
f(1\sw 1)\cdot c\su 2 = f(c),
$$
$$
\theta^{-1}\circ\theta(m) = 1\sw 0\cdot \theta(m)(1\sw 1) = 1\sw 01\sw 1\su
1\cdot m\cdot 1\sw 1\su 2 = m\cdot 1 = m.
$$

(2) For any $n>1$ consider the map $h^n: C^n_\psi(A,M) \to
C^{n-1}_\psi(A,M)$, given by 
$$
h^n(f)(c,a^1,\ldots, a^{n-1}) = c\su 11\sw 0\cdot f(1\sw 1,c\su 2,
a^1,\ldots, a^{n-1}).
$$
We will show that $h$ is a contracting homotopy, i.e., $h^{n+1}
d^{n} + d^{n-1}h^n = C_\psi^n(A,M)$. We use properties of the
translation map listed above and the definitions of $d$ and $h$ to
compute
\begin{eqnarray*}
d^{n-1}h^n (f)(c,a^1,\ldots, a^n) & = & c\su 1(c\su 2 a^1)\sw 0\cdot
h^n (f)((c\su 2 a^1)\sw 1, a^2,\ldots, a^n)\\
&&\hspace{-4cm} + \sum_{k=1}^{n-1}(-1)^k h^n(f)(c,
a^1,\ldots, a^ka^{k+1},\ldots, a^n) +(-1)^n
h^n(f)(c,a^1,\ldots,a^{n-1})\cdot a^n\\
&&\hspace{-4cm} = c\su 1(c\su 2 a^1)\sw 0(c\su 2 a^1)\sw 1\su 1 
1\sw 0\cdot
f(1\sw 1,(c\su 2 a^1)\sw 1\su 2, a^2,\ldots, a^n)\\
&&\hspace{-4cm} + \sum_{k=1}^{n-1}(-1)^k c\su 11\sw 0\cdot f(1\sw 1,c\su 2,
a^1,\ldots, a^ka^{k+1},\ldots, a^n)\\
&&\hspace{-4cm}  +(-1)^n c\su 11\sw 0\cdot f(1\sw 1,
c\su 2,a^1,\ldots,a^{n-1})\cdot a^n\\
&&\hspace{-4cm} = c\su 1 1\sw 0\cdot
f(1\sw 1,c\su 2 a^1, a^2,\ldots, a^n)+(-1)^n c\su 11\sw 0\cdot f(1\sw 1,
c\su 2,a^1,\ldots,a^{n-1})\cdot a^n\\
&&\hspace{-4cm} +
\sum_{k=1}^{n-1}(-1)^k c\su 11\sw 0\cdot f(1\sw 1,c\su 2,
a^1,\ldots, a^ka^{k+1},\ldots, a^n)\\
 &&\hspace{-4cm} = - c\su 11\sw 0\cdot d^nf(1\sw 1,c\su 2, a^1,\ldots,
a^n) + c\su 11\sw 0c\su 2_\alpha\cdot f(1\sw 1^\alpha, a^1,\ldots,
a^n)\\
 &&\hspace{-4cm} = - h^{n+1}(d^nf)(c, a^1,\ldots,
a^n) + c\su 1c\su 2\sw 0\cdot f(c\su 2\sw 1, a^1,\ldots,
a^n)\\
&&\hspace{-4cm} = - h^{n+1}(d^n f)(c, a^1,\ldots,
a^n) + f(c, a^1,\ldots,
a^n).
 \end{eqnarray*}
Therefore $h^\bullet$ is a contracting homotopy, so that for all $n>
0$, $H^n_\psi(A,M) =0$ as claimed.
\end{proof}

In particular, Proposition~\ref{prop.Galois} implies that if $A$ is a
Hopf algebra and $\psi:A\otimes A\to A\otimes A$, $a\otimes a'\to a'\sw
1\otimes aa'\sw 2$, then $H^0_\psi(A,M) = M$ and $H^n_\psi(A,M) = 0$ for
all $n>0$.   Furthermore, in view of 
Proposition~\ref{prop.proj},
Proposition~\ref{prop.Galois} implies that for a $C$-Galois object $A$,
$A\otimes C$ is a projective $A$-bimodule. Notice also that
Proposition~\ref{prop.Galois}(1) states that the entwined zero-cohomology
group of $A(B)^C$ is the same as the Hochschild zero-cohomology
group of $B$.

As another example we compute the zero cohomology group of
an $A$-Galois coextension with values in $A$.

\begin{example}
Suppose $(A,C)_\psi$ is the canonical entwining structure associated to
an $A$-Galois coextension
$C(B)_A$. 
Then $H^0_\psi(A,A) = {}^B{\rm End}_A(C)$ (the space of left $B$-comodule
right $A$-module endomorphisms of $C$).
\end{example}
\begin{proof}
In general 
$H^0_\psi(A,A) = \{\phi\in {\rm Hom}(C,A)\; |
\; \forall c\in C, a\in A \;\; a_\alpha\phi(c^\alpha) =
\phi(c)a\}. 
$
Now, \cite[Theorem~2.4*]{Brz:mod} yields the assertion. 
\end{proof}

\begin{sloppy}
\section{Comodule valued cohomology of  an 
entwining structure}
\end{sloppy}
The construction of the previous section can be dualised to produce a
$\psi$-twisted cohomology of a coalgebra. Thus the aim of this section
is to describe a cochain complex associated to an entwining structure
$(A,C)_\psi$ and a $C$-bicomodule $V$.

Consider the cobar resolution of a coalgebra $C$,
$\C(C) = (\C^\bullet(C),\bar{\delta}) $, where 
$$
\C^n(C) = C^{n+2},\quad \bar{\delta}^n =
\sum_{k=0}^{n} (-1)^k C^{k}\otimes\Delta\otimes C^{n-k}:
C^{n+2}\mapsto C^{n+3}.
$$
One easily checks that $\bar{\delta}$ is a $C$-bicomodule map so
that one can define the
cochain complex $\C_\psi(C) = (\C_\psi^\bullet(C),\bar{\delta})$ via
$ \C_\psi(C) = (C\otimes A)\square_C
\C(C)$, where $\square_C$ denotes the cotensor product over $C$.
Explicitly, $\C_\psi^n(C) = C\otimes A\otimes C^{n+1}$ and 
$$\bar{\delta}^n =
(C\otimes \psi\otimes C^n)\circ(\Delta\otimes A\otimes C^n)
+\sum_{k=1}^n(-1)^k C\otimes A\otimes C^{k-1}\otimes\Delta\otimes
C^{n-k}.
$$
As in the case of the $\psi$-twisted bar resolution $\B^\psi(A)$ one
has,
\begin{lemma}
$\C_\psi(C)$  is a resolution of $C$ and $C\otimes A$. Furthermore, $\bar{\delta}$ is a
$C$-bicomodule map.
\end{lemma}
\begin{proof} The contracting homotopy is $h^n: \C_\psi^n(C)\to
\C_\psi^{n-1}(C)$, $h^n = (-1)^{n+1}C\otimes A\otimes C^{n}\otimes \eps$,
while the augmentations are $(C\otimes 1\otimes C)\circ\Delta$ and
$(C\otimes\psi)\circ(\Delta\otimes A)$.
\end{proof}

Let $V$ be a $C$-bicomodule. 
Define the cochain complex $A_\psi (C,V) = 
(A_\psi(C,V)^\bullet, \bar{d})$ by $ A_\psi (C,V) ={}^C{\rm
Hom}^C(V,\C_\psi(C))$, where ${}^C{\rm Hom}^C$ denotes the Hom-bifunctor
from the category of 
$C$-bicomodules to the category of $k$-vector spaces.
Using the natural identification
$
^C{\rm
Hom}^C(C,C\otimes A\otimes C^{n+1}) = {\rm Hom}(V,A\otimes C^n),
$
one explicitly obtains
$$
A_\psi^n(C,V) = {\rm Hom}(V,A\otimes C^n), \quad {\bar d}^n 
: {\rm Hom}(V,A\otimes
C^n) \to  {\rm Hom}(V,A\otimes C^{n+1}),
$$
$$ \bar{d}^n f =
(C\otimes f)\circ(\psi\otimes C^n)\circ {}^V\rho +
\sum_{k=1}^n(-1)^k 
(A\otimes C^{k-1}\otimes\Delta\otimes C^{n-k})\circ f
+(-1)^{n+1}(f\otimes C)\circ\rho^V,
$$
where $\rho^V$ (${}^V\rho$ resp.) denotes the right (left) coaction of $C$
on $V$.

The Cartier complex over $C$ with values in
$V$ (cf.\ \cite{Car:coh}) 
is included in the complex $A_\psi(C,V)$. More
precisely one has
\begin{lemma}
Let $C(C,V)$ be the Cartier complex over
$C$ with values in $V$. 
Then the map $\bar{j}:C(C, V) \to A_\psi(C, V)$ given
by 
$$
\bar{j}^n :{\rm Hom}(V,C^n)\to {\rm Hom}(V,A\otimes C^n), \qquad \bar{j}^n:
f\mapsto 1\otimes f,
$$
is a monomorphism of differential complexes.
\label{lemma.car}
\end{lemma}

The cohomology of the complex $A_\psi(C, V)$ is
denoted by $H_\psi(C, V)$ and is called an {\em entwined
cohomology of $C$ with values in $V$}.

As an example we compute the (canonical) entwined cohomology of an 
$A$-Galois coextension.
\begin{proposition}
Let $(A,C)_\psi$ be the canonical entwining structure associated to an
algebra-Galois coextension $C(B)_A$ in Example~\ref{cocan.ex}, and let $V$
be a $C$-bicomodule. Then

(1)  $H_\psi^0(C,V) = V_B := \{v\in V\; | \;  (\pi\otimes
V)\circ{}_V\rho(v) = (V\otimes\pi)\circ\rho_V(v)\}$. Here $\pi:C\to
B=C/I$ is the canonical epimorphism.

(2) If $B=k$, then $H_\psi^n(C,V) =
0$, for all $n>0$.
\label{prop.coGalois}
\end{proposition}
\begin{proof}
The proof is dual to that of Proposition~\ref{prop.Galois}. The
isomorphism in part (1) is
$$
H_\psi^0(C,V)\to V_B,\quad f\mapsto \eps\circ\rho_C\circ(C\otimes
f)\circ {}^V\!\!\rho,
$$
while the contracting homotopy in part
(2) is $h^n: A_\psi^n(C,V) \to
A_\psi^{n-1}(C,V)$,  
$$
h^n(f) = (\hat{\tau}\otimes
C^{n-1})(C\otimes\eps\otimes C^{n})(C\otimes\rho_C\otimes
C^n)(\Delta\otimes f){}^V\!\!\rho,
$$
where $\hat{\tau} = (\eps\otimes A)\circ{\rm cocan}^{-1}$ is the {\em
cotranslation map}.
\end{proof}

Yet another example is the zero-cohomology group of the canonical entwining
structure associated to a $C$-Galois extension.
\begin{example}
Suppose $(A,C)_\psi$ is the canonical entwining structure associated to
a $C$-Galois extension
$A(B)^C$. View $C$ as a $C$-bicomodule via the coproduct. 
Then $H_\psi^0(C,C) = {}_B{\rm End}^C(A)$ (the space of left $B$-module
right $C$-comodule endomorphisms of $A$).
\end{example}
\begin{proof}
In general 
$H_\psi^0(C,C)= \{\phi\in {\rm Hom}(C,A)\; |
\; \forall c\in C, \; \phi(c\sw 2)_\alpha\otimes c\sw 1^\alpha =
\phi(c\sw 1)\otimes c\sw 2\}. 
$
Now, \cite[Theorem~2.4]{Brz:mod} yields the assertion. 
\end{proof}

 \section{Cup products on $C_\psi(A)$.} 
The structure of the entwined cohomology of an algebra is particularly
rich when the cohomology takes its values in the algebra itself. The
situation is very much reminiscent of the Hochschild cohomology of an
algebra, as described in \cite{Ger:coh}. 
In this
section we study the structure of $C_\psi(A):= C_\psi(A, A)$ along
the lines of \cite{Ger:coh}. 

There are (at least) two ways of defining an associative algebra
structure or cup products on  $C_\psi(A)$.
For any $f\in C_\psi^m(A)$, $g\in C_\psi^n(A)$ define
$$
f\cup g =\mu\circ(f\otimes g)\circ (\rho^m_R\otimes A^n) \in
C_\psi^{m+n}(A),
$$
where $\rho^m_R$ is the right coaction of $C$ on $C\otimes A^m$
described in the introduction. Using the $\alpha$-notation for the
entwining map we can write explicitly,
$$
(f\cup g)(c,a^1,\ldots, a^{m+n}) = f(c\sw 1,
a^1_{\alpha_1},\ldots ,a^m_{\alpha_m})g(c\sw
2^{\alpha_1\ldots\alpha_m},a^{m+1},\ldots ,a^{m+n}).
$$
We have
\begin{lemma}
$C_\psi(A)$ is a graded
associative algebra with the product $\cup$. Furthermore, $d$ is a
degree 1 derivation in this algebra, i.e., for all
 $f\in C_\psi^m(A)$, $g\in C_\psi^n(A)$,
$$
d(f\cup g) = df\cup g + (-1)^mf\cup dg.
$$
\label{cup.lemma}
\end{lemma}
\begin{proof}
This follows from Proposition~\ref{prop.ascup} and
Proposition~\ref{c.comp.prop} below, but it can also be
 proven by straightforward manipulations with the definition of
an entwining structure. In particular, the associativity of $\cup$
follows from the right pentagon in the bow-tie diagram, while to prove
the derivation property of $d$ one needs to use the definition of $d$
and both pentagons in the bow-tie diagram.
\end{proof}

Lemma~\ref{cup.lemma} implies that the cup product $\cup$ defines 
the product in the cohomology $H_\psi(A):=H_\psi(A,A)$. This
product is also denoted by $\cup$.  

The second type of a cup product is defined as follows. For any $f\in 
C_\psi^m(A)$, $g\in
C_\psi^n(A)$, consider
$$
f\sqcup g= \mu\circ(A\otimes g)\circ (\psi\otimes A^n)\circ(C\otimes
f\otimes A^n)\circ(\Delta\otimes A^{m+n})\in C_\psi(A)^{m+n}.
$$
Explicitly,
$$
(f\sqcup g)(c,a^1,\ldots ,a^{m+n}) = f(c\sw 2,a^1,\ldots ,
a^m)_\alpha g(c\sw
1^\alpha, a^{m+1},\ldots ,a^{m+n}).
$$
\begin{lemma}
$C_\psi(A)$ is an
associative algebra with the product $\sqcup$. 
Furthermore, $d$ is a
degree 1 derivation in this algebra.
\label{sqcup.lemma}
\end{lemma}
\begin{proof}
This lemma follows from Proposition~\ref{prop.ascup} and
Proposition~\ref{c.comp.prop} below too, but it can also be proven by
 a straightforward application of definitions of $d$ and
$\psi$, which uses both pentagons of the bow-tie diagram.
\end{proof}

Lemma~\ref{sqcup.lemma} implies that the product $\sqcup$ 
defines the product in $H_\psi(A)$. This
product is also denoted by $\sqcup$. 

Both cup products are closely related to each other 
in the cohomology $H_\psi(A)$. 
This is described in the following
\begin{theorem}
For all
cohomology classes
 $\xi\in H_\psi^m(A)$, $\eta\in H_\psi^n(A)$, 
$$
\xi\cup\eta = (-1)^{mn}\eta\sqcup\xi.
$$
\label{cup.theorem}
\end{theorem}

Theorem~\ref{cup.theorem} is a generalisation of the result of
Gerstenhaber for the Hochschild cohomology of an algebra in
\cite{Ger:coh}. Indeed, take the trivial entwining structure
$(A,k)_\psi$ (i.e., $\psi$ is a usual flip canonically identified with
the identity automorphism of $A$). In this case $C_\psi(A)$ is simply the
Hochschild complex of $A$ and both the cup products $\cup$ and $\sqcup$
become the standard cup product in the Hochschild complex.
Theorem~\ref{cup.theorem} is then a simple consequence of \cite[Section~7,
Corollary~1]{Ger:coh}. 

The rest of this section is devoted to the proof of
Theorem~\ref{cup.theorem}. We employ a method similar to the one used in
\cite{Ger:coh}. The method used there is based on pre-Lie systems and
pre-Lie algebras. Although we are not able to associate a pre-Lie system
to the general complex $C_\psi(A)$, still it is possible to construct a
system whose properties suffice to prove  Theorem~\ref{cup.theorem}.
First we introduce the following generalisation of the notion of a comp
algebra introduced in \cite{GerSch:alg}\cite{Ger:coh} 
\begin{definition}
A {\em (right)
weak comp algebra} $(V^\bullet,\diamond,\pi)$  consists of a sequence of vector
spaces $V^0, V^1,V^2,\ldots$, an element $\pi\in V^2$ and $k$-linear
operations $\diamond_i : V^m\otimes V^n\to V^{m+n-1}$ for $i\geq 0$ such
that for any $f\in V^m$, $g\in V^n$, $h\in V^p$,

(1) $f\diamond_i g =0$ if $i>m-1$;

(2) $(f\diamond_ig)\diamond_jh =f\diamond_i (g\diamond_{j-i} h)$  
if $i\leq j <n+i$;

(3) if either $g=\pi$ or $h=\pi$,
$$
(f\diamond_ig)\diamond_jh = 
  (f\diamond_j h)\diamond_{i+p-1} g \quad \mbox{if $j <i$};
$$

(4) $\pi\diamond_0\pi = \pi\diamond_1\pi$.
\label{comp.def}
\end{definition}
The weak comp algebra is a (right) comp algebra in the sense of
Gerstenhaber and Schack \cite{GerSch:alg} if the condition (3) 
holds for all $g$ and $h$. In the case of a comp algebra 
one can define an associative
cup product in $V=\bigoplus_{i=0}V^i$. We shall see that this cup
product splits into two cup products for the weak comp algebra. 
Notice that condition (3) of Definition~\ref{comp.def}
 implies also that if either $g=\pi$ or $h=\pi$,
$$
(f\diamond_ig)\diamond_jh = 
  (f\diamond_{j-n+1} h)\diamond_{i} g \quad \mbox{if $j\geq n+i$}.
$$

Given a weak comp algebra $(V^\bullet,\diamond,\pi)$ define an operation
$\diamond: V^m\otimes V^n \to V^{m+n-1}$, for all $f\in
V^m$, $g\in V^n$ given  by
$$
f\diamond g = \sum_{i=0}^{m-1}(-1)^{i(n-1)} f\diamond _i g.
$$

Furthermore one can introduce two additional operations $\cup,\sqcup:
V^m\otimes V^n \to V^{m+n}$ for all $f\in
V^m$, $g\in V^n$ given by
\begin{equation}
f\cup g = (\pi\diamond_0f)\diamond_m g, \qquad 
f\sqcup g = (\pi\diamond_1g)\diamond_0f .
\label{cups}
\end{equation}
The importance of these operations is revealed in the following
\begin{proposition}
Given a weak comp algebra $(V^\bullet,\diamond,\pi)$,
$V=\bigoplus_{i=0}V^i$  is a (non-unital) graded associative algebra with respect to
each of the operations $\cup$ and $\sqcup$ defined in equations (\ref{cups}). 
\label{prop.ascup}
\end{proposition}
\begin{proof}
We prove this proposition for the operation $\cup$, the proof for
$\sqcup$ is analogous. Take any $f\in V^m$, $g\in V^n$ and $h\in V^p$. To
prove the associativity of $\cup$ we need to show that $(f\cup g)\cup h
= (\pi\diamond_0((\pi\diamond_0f)\diamond_mg))\diamond_{m+n} h$ equals to
$f\cup(g\cup h) = (\pi\diamond_0 f)\diamond_m((\pi\diamond_0
g)\diamond_n h)$. First using condition (2) in Definition~\ref{comp.def}
one easily finds that $(\pi\diamond_0 f)\diamond_m((\pi\diamond_0
g)\diamond_n h) = ((\pi\diamond_0 f)\diamond_m (\pi\diamond_0
g))\diamond_{m+n} h$, so that only the equality
$\pi\diamond_0((\pi\diamond_0f)\diamond_mg) = (\pi\diamond_0 f)
\diamond_m(\pi\diamond_0
g)$ needs to be shown. We have
\begin{eqnarray*}
\pi\diamond_0((\pi\diamond_0f)\diamond_mg) &=&
(\pi\diamond_0(\pi\diamond_ 0 f))\diamond_m g
= ((\pi\diamond_0 \pi)\diamond_ 0 f)\diamond_m g
= ((\pi\diamond_1 \pi)\diamond_ 0 f)\diamond_m g\\
&=& ((\pi\diamond_0 f)\diamond_m \pi)\diamond_m g
= (\pi\diamond_0 f)\diamond_m(\pi\diamond_0 g),
\end{eqnarray*}
where we used Definition~\ref{comp.def}(2) to derive the first, second
and the fifth equalities, Definition~\ref{comp.def}(4) to derive the
third equality and Definition~\ref{comp.def}(3) to obtain the fourth
one.
\end{proof}

Conditions in Definition~\ref{comp.def}
allow one to use the same method as in the proof of 
\cite[Theorem~2]{Ger:coh}
to prove the following
\begin{proposition}
Let $(V^\bullet,\diamond,\pi)$ be a weak comp algebra. Then 
for all $f\in V^m$, $g\in V^n$ and $h\in V^p$ 
we have:

(1) if either $f$, $g$ or $h$ is equal to $\pi$ then
$$
(f\diamond g)\diamond h - f\diamond(g\diamond h) =
\sum (-1)^{i(n-1)+j(p-1)}(f\diamond_i g)\diamond_j h,
$$
 where the sum is over those $i$ and $j$ with either $0\leq j\leq i-1$ or
$n+i\leq j\leq m+n-2$.

(2) $(f\diamond g)\diamond \pi - f\diamond(g\diamond \pi) = (-1)^{n-1}
\left( (f\diamond \pi)\diamond g - f\diamond(\pi\diamond g)\right).
$
\label{pre-Lie.prop}
\end{proposition}
Proposition~\ref{pre-Lie.prop} allows one to construct a coboundary in a
weak comp algebra.
\begin{proposition}
A weak comp algebra $(V^\bullet,\diamond,\pi)$ is a cochain complex with
a coboundary $d:V^m\to V^{m+1}$,
$
df = (-1)^{m-1}\pi\diamond f -f\diamond\pi.
$
Furthermore, $d$ is a degree one derivation in both algebras 
$(V,\cup)$ and $(V,\sqcup)$
\label{prop.inner}
\end{proposition}
\begin{proof}
The first part is a simple consequence of Proposition~\ref{pre-Lie.prop}
and  easily verifiable fact that $\pi\diamond\pi = 0$. The second
part can be proven by direct computation. We display it for the cup
product $\cup$. Explicitly, for any $f\in V^m$, $g\in V^n$ 
 one needs to show that
$
d(f\cup g) = df\cup g +(-1)^mf\cup dg.
$
Using definitions of $\cup$, $\diamond$ and $d$ this amounts to showing that
$$
\Gamma_1 = 
(-1)^{m+n-1}\pi\diamond_0((\pi\diamond_0f)\diamond_m g)
+\pi\diamond_1((\pi\diamond_ 0f)\diamond_mg) - \sum_{j=0}^{m+n-1} (-1)^j
((\pi\diamond_0f)\diamond_m g)\diamond_j\pi
$$
is equal to
\begin{eqnarray*}
\Gamma_2 &=& (-1)^{m+n-1} (\pi\diamond_0f)\diamond_m(\pi\diamond_0g)
+(\pi\diamond_0(\pi\diamond_1 f))\diamond _{m+1} g \\
&&-\sum_{i=0}^{m-1}(-1)^i (\pi\diamond_0(f\diamond_i\pi))\diamond_{m+1} g
-(-1)^m\sum_{k=0}^{n-1}(-1)^k(\pi\diamond_0 f)\diamond_m(g\diamond_k\pi)\\
&& +
(-1)^{m-1}( (\pi\diamond_0(\pi\diamond_0f))\diamond_{m+1}g -
(\pi\diamond_0 f)\diamond_m (\pi\diamond_1 g)). 
\end{eqnarray*}
The first term in $\Gamma_1$ equals the first term in $\Gamma_2$ by the
same calculation as in the proof of Proposition~\ref{prop.ascup}. Using
a chain of arguments as in the proof of
Proposition~\ref{prop.ascup} but without the fourth step, one easily shows that the second term in
$\Gamma_1$ is the same as the second term in $\Gamma_2$. Again, a part
of the argument in the proof of Proposition~\ref{prop.ascup} allows one
to transform the first term inside the final brackets in $\Gamma_2$ to
the form $((\pi\diamond_0 f)\diamond_m\pi)\diamond_{m+1}g$ and then
Definition~\ref{comp.def}(2) implies that the last line in $\Gamma_2$
vanishes. Now consider the third term in $\Gamma_1$. If $j\leq m-1$ then
conditions (3) and (2) in Definition~\ref{comp.def} imply that
$$
((\pi\diamond_0 f)\diamond_m g)\diamond_j\pi =
((\pi\diamond_0f)\diamond_j\pi)\diamond_{m+1} g = (\pi\diamond_0(f\diamond_j\pi))\diamond_{m+1}
g ,
$$
so that the part of the sum in the last term in $\Gamma_1$ for $j\leq
m-1$ is the same as the third term in $\Gamma_2$. Similarly,
Definition~\ref{comp.def}(2) implies that the remaining part of this sum
is the same as the fourth term in $\Gamma_2$. This completes the proof
that $d$ is a derivation in the algebra $(V,\cup)$. Similar
arguments show that $d$ is a derivation in the algebra
$(V,\sqcup)$. 
\end{proof}

The relationship between three operations $\diamond$, $\cup$ and
$\sqcup$ is revealed in the following
\begin{theorem}
Let $(V^\bullet,\diamond,\pi)$ be a weak comp algebra. Then 
for all $f\in V^m$, $g\in V^n$:
$$
f\diamond dg - d(f\diamond g) + (-1)^{n-1} df\diamond g =
(-1)^{n-1}\left( g\sqcup f - (-1)^{mn} f\cup g\right).
$$
\label{pre-Lie.thm}
\end{theorem}
\begin{proof}
Expand the left hand side of the above equality using definition of $d$
in Proposition~\ref{prop.inner}.
This produces six terms, four of which cancel because of
Proposition~\ref{pre-Lie.prop}(2). One is left to show that
$
(\pi\diamond f)\diamond g - \pi\diamond(f\diamond g) = 
(-1)^{m-1} \left( g\sqcup f - (-1)^{mn} f\cup
g\right) .
$
By Proposition~\ref{pre-Lie.prop}(1), 
$
(\pi\diamond f)\diamond g - \pi\diamond(f\diamond g) = 
(-1)^{m-1}((\pi\diamond_1 f)\diamond_0 g -
(-1)^{mn}(\pi\diamond_0f)\diamond_m g),
$
therefore the required equality holds once the definition of cup
products in equation~(\ref{cups}) is taken into account. 
\end{proof}

By Proposition~\ref{prop.inner} the cohomology of a
weak comp algebra can be equipped with two algebra structures 
corresponding to
products $\cup$ and $\sqcup$. Now Theorem~\ref{pre-Lie.thm} implies

\begin{corollary}
Let $(V^\bullet,\diamond,\pi)$ be a weak comp algebra and let $H(V)$ 
denote its cohomology with respect to the coboundary operator defined 
in Proposition~\ref{prop.inner}. Then for all
cohomology classes
 $\xi\in H^m(V)$, $\eta\in H^n(V)$, 
$$
\xi\cup\eta = (-1)^{mn}\eta\sqcup\xi.
$$
\label{cup.cor}
\end{corollary}
\begin{proof}
If $f$ and $g$ are cocycles one has 
$
d(f\diamond g)  =
(-1)^{n}\left( g\sqcup f - (-1)^{mn} f\cup g\right) ,
$
and the corollary follows.
\end{proof}

Now we are ready to prove Theorem~\ref{cup.theorem} by associating a 
weak comp algebra to the complex $C_\psi(A)$. 
First, for any $f\in C_\psi^m(A)$, $g\in C_\psi^n(A)$,
 define $f\diamond_i g \in  C_\psi^{m+n -1}(A)$ by
\begin{equation}
f\diamond_i g = \left\{ \begin{array}{ll} 
            f\circ (C\otimes A^i \otimes g\otimes A^{m-i-1})\circ
(\rho^i_R \otimes A^{m+n-i-1}) & \mbox{if $0\leq i <m$}\\
0 & \mbox{otherwise}
\end{array}
\right. 
\label{diamond}
\end{equation}
Explicitly for all $0\leq i<m$,
\begin{eqnarray*}
&&(f\diamond_i g)(c,a^1,\ldots,a^{m+n-1}) \\&&{ }\qquad\qquad =
f(c\sw 1,a^1_{\alpha_1}, \ldots, a^i_{\alpha_i},
g(c\sw 2^{\alpha_1\ldots \alpha_i}, a^{i+1},\ldots,a^{i+n}),
a^{n+i+1},
\ldots, a^{m+n-1}).
\end{eqnarray*}
Next consider the two-coboundary $\pi\in
C_\psi^2(A)$ given by $\pi = \eps\otimes\mu$. It is a
coboundary, since one easily checks that 
$\pi = d (\eps\otimes A)$. 
\begin{proposition}
$(C_\psi^\bullet(A),\diamond,\pi)$ is a weak comp algebra.
\label{c.comp.prop}
\end{proposition}
\begin{proof}
Condition (1) in Definition~\ref{comp.def} is clearly satisfied. 
Definition~\ref{comp.def}(2) can be proven by a straightforward
calculation which uses the right pentagon in the bow-tie diagram. 
Next, notice that for any 
$f\in C_\psi^m(A)$
$$
f\diamond_i\pi= f\circ(C\otimes A^i\otimes\mu\otimes A^{m-i-1}),
$$
take $0\leq j\leq i-1$ and compute
 \begin{eqnarray*}
(f\diamond_ig)\diamond_j\pi & = & (f\diamond_ig)(C\otimes A^j\otimes
\mu\otimes A^{m+n-j-1})\\
&& \hspace{-3cm} =  f(C\otimes A^i\otimes g\otimes A^{m-i-1})
(\rho_R^i\otimes
A^{m+n-i-1})(C\otimes A^j\otimes
\mu\otimes A^{m+n-j-1})\\
&& \hspace{-3cm} = f(C\otimes A^i\otimes g\otimes A^{m-i-1})
(C\otimes A^j\otimes\mu
\otimes A^{i-j-1}\otimes C\otimes A^{m+n-i})
(\rho_R^{i+1}\otimes A^{m+n-i})\\
&& \hspace{-3cm} = f(C\otimes A^j\otimes \mu\otimes A^{m-j-1})
(C\otimes A^{i+1}
\otimes g\otimes A^{m-i-1})(\rho_R^{i+1}\otimes A^{m+n-i})\\
&& \hspace{-3cm} =  (f\diamond_j\pi)(C\otimes A^{i+1}
\otimes g\otimes A^{m-i-1})(\rho_R^{i+1}\otimes A^{m+n-i})\\
&& \hspace{-3cm} = (f\diamond_j\pi)\diamond_{i+1}g,
\end{eqnarray*}
where we used Lemma~\ref{lemma.system} to derive the third equality.
This proves Definition~\ref{comp.def}(3) with $h=\pi$. The case $g=\pi$
is proven in a similar way. The proof of condition
Definition~\ref{comp.def}(4) is again straightforward.
\end{proof}

One can easily verify that the cup products in $C_\psi(A)$ defined at
the beginning of this section, are given by
equations (\ref{cups}) while the coboundary operator is given  by
the formula in Proposition~\ref{prop.inner}, with operations
$\diamond_i$ defined in equation (\ref{diamond}). 
Therefore Theorem~\ref{cup.theorem} immediately follows from 
Corollary~\ref{cup.cor}.

Dually, one can associate a weak comp algebra to a $C$-valued entwined 
cochain complex  of $C$, $A_\psi(C) = A_\psi(C,C)$. In this case the
operations $\diamond_i$ are defined as
$$
f\diamond_i g = \left\{ \begin{array}{ll} 
            (\rho_i^R \otimes C^{m+n-i-1})
\circ(A\otimes C^i \otimes g\otimes C^{m-i-1})\circ f  
& \mbox{if $0\leq i <m$}\\
0 & \mbox{otherwise}
\end{array}
\right. 
$$
for all $f\in A_\psi^m(C)$, $g\in A_\psi^n(C)$, while $\pi =
1\otimes\Delta$. 

\section{The $\psi$-equivariant cohomology of $A$}
Guided by the results of the previous section, we construct here a
subcomplex $C_{\psi-e}(A)$ of $C_\psi(A)$ whose cohomology has a
graded-commutative algebra structure given by the cup product.

For any $n\in {\bf Z}_{\geq 0}$ consider a vector subspace of
$C_\psi^n(A)$, 
$$
C_{\psi-e}^n(A) = \{f\in {\rm Hom}(C\otimes A^n,A)\; |\; (f\otimes
C)\circ\rho^n_R = \psi\circ(C\otimes f)\circ\rho^n_L\}\subseteq 
C_\psi^n(A).
$$
\begin{lemma}
$(C_{\psi-e}^\bullet(A), \diamond,\pi)$ is a weak comp algebra with
$\diamond_i$ given by equation (\ref{diamond}) and $\pi =
\eps\otimes\mu$.
\end{lemma}
\begin{proof}
The condition for $f$ to be in $C_{\psi-e}^n(A) $ can be explicitly 
written for all
$c\in C$, $a^1,\ldots, a^n\in A$ 
\begin{equation}
f(c\sw 1, a^1_{\alpha_1},\ldots, a^n_{\alpha_n})\otimes c\sw
2^{\alpha_1\cdots\alpha_n} = f(c\sw 2,a^1,\ldots,a^n)_\alpha\otimes c\sw
1^\alpha.
\label{finbar}
\end{equation}
If $f=\pi$ this is precisely the left pentagon in the bow-tie diagram.
Therefore $\pi\in C_{\psi-e}^2(A)$. We will show that for all $f\in
C_{\psi-e}^m(A)$, $g\in C_{\psi-e}^n(A) $, $f\diamond_ig
\in C_{\psi-e}^{m+n-1}(A)$. Take any  $c\in C$, 
$a^1,\ldots, a^{m+n-1}\in A$ and compute
\begin{eqnarray*}
f\diamond_i g(c\sw 1, a^1_{\alpha_1},\ldots, 
a^{m+n-1}_{\alpha_{m+n-1}})\otimes c\sw
2^{\alpha_1\cdots\alpha_{m+n-1}} &=& \\
&& \hspace{-7.5cm} = f(c\sw 1, a^1_{\alpha_1\beta_1},\ldots,
a^i_{\alpha_i\beta_i}, g(c\sw 2^{\beta_1\ldots\beta_i},
a^{i+1}_{\alpha_{i+1}},\ldots
,a^{n+i}_{\alpha_{n+i}}),a^{n+i+1}_{\alpha_{n+i+1}},
\ldots , a^{m+n-1}_{\alpha_{m+n-1}})\\
&&\otimes c\sw
3^{\alpha_1\cdots\alpha_{m+n-1}}\\
&& \hspace{-7.5cm} = f(c\sw 1, a^1_{\alpha_1},\ldots,
a^i_{\alpha_i}, g(c\sw 2^{\alpha_1\ldots\alpha_i}\sw 1,
a^{i+1}_{\alpha_{i+1}},\ldots
,a^{n+i}_{\alpha_{n+i}}),a^{n+i+1}_{\alpha_{n+i+1}},\ldots 
,a^{m+n-1}_{\alpha_{m+n-1}})\\
&&\otimes c\sw
2^{\alpha_1\ldots\alpha_i}\sw 2^{\alpha_{i+1}\ldots \alpha_{m+n-1}}\\
&& \hspace{-7.5cm} = f(c\sw 1, a^1_{\alpha_1},\ldots,
a^i_{\alpha_i}, g(c\sw 2^{\alpha_1\ldots\alpha_i}\sw 2,
a^{i+1},\ldots ,a^{n+i})_\beta,a^{n+i+1}_{\alpha_{n+i+1}},\ldots 
,a^{m+n-1}_{\alpha_{m+n-1}})\\
&&\otimes c\sw 2^{\alpha_1\ldots\alpha_i}
\sw 1^{\beta\alpha_{n+i+1}\ldots\alpha_{m+n-1}}\\
&& \hspace{-7.5cm} = f(c\sw 1, a^1_{\alpha_1\beta_1},\ldots,
a^i_{\alpha_i\beta_i}, g(c\sw 3^{\alpha_1\ldots\alpha_i},
a^{i+1},\ldots ,a^{n+i})_\beta,a^{n+i+1}_{\alpha_{n+i+1}},\ldots 
,a^{m+n-1}_{\alpha_{m+n-1}})\\
&&\otimes c\sw 2^{\beta_1\ldots\beta_i\beta
\alpha_{n+i+1}\ldots\alpha_{m+n-1}}\\
&& \hspace{-7.5cm} = f(c\sw 2, a^1_{\alpha_1},\ldots,
a^i_{\alpha_i}, g(c\sw 3^{\alpha_1\ldots\alpha_i},
a^{i+1},\ldots ,a^{n+i}),a^{n+i+1},\ldots 
,a^{m+n-1})_\alpha\otimes c\sw 1^\alpha\\
&& \hspace{-7.5cm} = f\diamond_i g(c\sw 1, a^1,\ldots, 
a^{m+n-1})_\alpha\otimes c\sw 2^\alpha ,
\end{eqnarray*}
where the right pentagon in the bow-tie diagram has been used in
derivation of the second and fourth equalities, and equation
(\ref{finbar}) for $f$ and $g$ in derivation of the third and fifth
equalities. 
\end{proof}

Therefore $(C_{\psi-e}^\bullet(A),\diamond,\pi)$ is a 
weak comp subalgebra
of $(C_{\psi}^\bullet(A),\diamond,\pi)$. Consequently, 
$C_{\psi-e}(A) = (C_{\psi-e}^\bullet(A),d)$ is a cochain subcomplex of
$C_\psi(A)$. The corresponding
cohomology is denoted by $H_{\psi-e}(A)$ and called {\em 
$\psi$-equivariant cohomology of $A$}. As for any weak comp algebra one
can define the cup products in $C_{\psi-e}(A)$ which will descend to the
cohomology $H_{\psi-e}(A)$. Notice, however, that 
$\cup$ coincides with $\sqcup$ in $C_{\psi-e}(A)$, and,
consequently in $H_{\psi-e}(A)$. Therefore, as a simple 
consequence of Theorem~\ref{cup.theorem} and
Theorem~\ref{pre-Lie.thm} one obtains the following 
\begin{theorem}
For all $f\in C_{\psi-e}^m(A)$, $g\in C_{\psi-e}^n(A)$:
$$
f\diamond dg - d(f\diamond g) + (-1)^{n-1} df\diamond g =
(-1)^{n-1}\left( g\cup f - (-1)^{mn} f\cup g\right) .
$$
Consequently the algebra $(H_{\psi-e}(A),\cup)$ is graded-commutative.
\end{theorem}

As an example of a $\psi$-equivariant complex we consider such a complex
associated to the canonical entwining structure.

Take any $A$-bimodule $M$ with the left and right actions ${}_M\rho$
and $\rho_M$ respectively, and consider two operations
$$
{\rm Hom}(C\otimes A^n,A)\otimes {\rm Hom}(C,M)\to {\rm Hom}(C\otimes
A^n,M), \quad f\otimes\phi \mapsto
f\cup\phi:={}_M\rho\circ(f\otimes\phi)\circ \rho^n_R,
$$
$$
 {\rm Hom}(C,M)\otimes {\rm Hom}(C\otimes A^n,A)\to {\rm Hom}(C\otimes
A^n,M), \quad \phi\otimes f \mapsto
\phi*f:=\rho_M\circ(\phi\otimes f)\circ \rho^n_L.
$$
\begin{example}
Let $(A,C)_\psi$ be the canonical entwining structure associated to a
$C$-Galois extension $A(B)^C$, and let $\tau = {\rm
can}^{-1}\circ(1\otimes C)$ be the translation map. 
Then $f\in C_{\psi-e}^n(A)$ if and only if $f\cup\tau = \tau*f$.
\end{example}
\begin{proof}
Apply ${\rm can}^{-1}$ to (\ref{finbar}) and then use the definitions of
the canonical entwining structure and the translation map to obtain
\begin{eqnarray*}
(f\cup\tau)(c,a^1,\ldots, a^n) &=& f(c\sw 2,a^1,\ldots,
a^n)_\alpha\tau(c\sw 1^\alpha)\\
&& \hspace{-3cm} = c\sw 1\su 1(c\sw 1\su 2f(c\sw 2,a^1,\ldots,
a^n))\sw 0\tau((c\sw 1\su 2f(c\sw 2,a^1,\ldots,
a^n))\sw 1)\\
&&\hspace{-3cm} =   \tau(c\sw 1)f(c\sw 2,a^1,\ldots,
a^n) = (\tau*f)(c,a^1,\ldots, a^n),
\end{eqnarray*}
where $\tau(c) = c\su 1\otimes c\su 2$.
\end{proof}

\section{Deformation  of entwining structures}
In this section we associate a double-complex to 
any entwining structure. The total cohomology of a
 particular subcomplex of this complex
gives the cohomological interpretation of deformation theory of
entwining structures.

Given an entwining structure $(A,C)_\psi$, let 
$C(A,C,\psi) = (C^{\bullet , \bullet}(A,C,\psi),d,\bar{d})$ 
be a double complex 
obtained by applying the Hom-functor in the category of $A$-bimodules,
$C$-bicomodules to the $\psi$-twisted bar and cobar resolutions,
i.e., 
$$
C(A,C,\psi) = {}_A^C{\rm Hom}_A^C(\B^\psi(A),\C_\psi(C)).
$$
Here $A\otimes C\otimes A^{n+1}$ is viewed as a $C$-bicomodule
via $(\sigma\otimes C\otimes A^{n+1})\circ(A\otimes \rho^n_L\otimes A)$
and $(A\otimes C\otimes A^n\otimes \sigma)\circ (A\otimes\rho^n_R\otimes
A)$, where $\sigma$ is a flip. Dually $C\otimes A\otimes C^{n+1}$ is
viewed as an $A$-bimodule via $(C\otimes \rho_n^L\otimes C)\circ
(\sigma\otimes A\otimes C^{n+1})$, and $(C\otimes \rho_n^R\otimes C)\circ
(C\otimes A\otimes C^{n+1}\otimes\sigma)$. Clearly bimodule structures
commute with the bicomodule structures. Explicitly we have
$$
C^{m,n}(A,C,\psi) = C_\psi^m(A,A\otimes C^n) = A_\psi^n(C,C\otimes A^m)
= {\rm Hom}(C\otimes A^m,A\otimes C^n),
$$
$$
d:C^{m,n}(A,C,\psi) \to C^{m+1,n}(A,C,\psi),\qquad \bar{d} : 
C^{m,n}(A,C,\psi)\to  C^{m+1,n}(A,C,\psi),
$$
where $d$ is a coboundary operator described in Section~2, corresponding
to the entwined cohomology of $A$ with values in the $A$-module
$(A\otimes C^n,\rho_n^L,\rho_n^R)$, while $\bar{d}$ is the coboundary
operator of
Section~3, corresponding to the entwined cohomology of $C$ with values in
the $C$-bicomodule $(C\otimes A^m, \rho^m_L,\rho^m_R)$. 
It follows directly from the construction  that 
$d\bar{d} = \bar{d}d$, so that $C(A,C,\psi)$ is a double cochain complex
as claimed, with the total coboundary operator
$D^{m,n} :
C^{m,n}(A,C,\psi) \to C^{m+1,n}(A,C,\psi)\oplus
C^{m,n+1}(A,C,\psi)$, $D^{m,n} = d +(-1)^m\bar{d}$. 

 From the point of view of deformation of entwining
structures the following modification of the double complex
$C(A,C,\psi)$
is of substantial importance.  Let
$
C_H^{m,0}(A, C,\psi)  =  {\rm Hom}(A^m,A)$, $ 
C_H^{0,n} (A,C,\psi)  =  {\rm Hom}(C,C^n)$, $
C_H^{m,n}(A,C,\psi)   =  {\rm Hom}(C^\otimes A^m,A\otimes C^n),
$ 
for
$m,n\geq 1$. In other words, the complex $C_H(A,C,\psi)$ is obtained
from  $C(A,C,\psi)$ by replacing the first line and first column
with the Hochschild complex of $A$ and the Cartier complex of 
$C$ viewed in $C(A,C,\psi)$  via the monomorphisms $j$ in
Lemma~\ref{lemma.hoch} and
$\bar{j}$ in Lemma~\ref{lemma.car}. 
Thus, explicitly, the complex $C_H(A,C,\psi)$ is

{\tiny
$$
\begin{diagram}
& & {\rm Hom}(A,A) & \rTo^{d} & {\rm Hom}(A^2,A) & \rTo^d 
&  {\rm Hom}(A^3,A) & \rTo^{d} \\
& & \dTo^{\bar{d}\circ j} & & \dTo^{\bar{d}\circ j}& &
\dTo^{\bar{d}\circ j}\\
{\rm Hom}(C,C)& \rTo^{d\circ \bar{j}} & {\rm Hom}(C\otimes A,A\otimes C) & \rTo^{d} & 
{\rm Hom}(C\otimes A^2, A\otimes C) & \rTo^{d} & {\rm Hom}
(C\otimes A^3, A\otimes C) & \rTo^{d} \\
\dTo^{\bar{d}}& & \dTo^{\bar{d}} & & \dTo^{\bar{d}}& & \dTo^{\bar{d}}\\
{\rm Hom}(C, C^2) & \rTo^{d\circ\bar{j}} & {\rm Hom}(C\otimes A, A\otimes C^2)& \rTo^{d}
& {\rm Hom}(C\otimes A^2, A\otimes C^2)& \rTo^{d} & {\rm Hom}
(C\otimes A^3, A\otimes C^2) & \rTo^{d} \\
\dTo^{\bar{d}}& & \dTo^{\bar{d}} & & \dTo^{\bar{d}}& & \dTo^{\bar{d}}\\
&&&&&&&&
\end{diagram}
$$}
This double complex combines into a complex
$C_H((A,C)_\psi) = (C_H^\bullet((A,C)_\psi),D)$,
$$
C_H^n((A,C)_\psi) = {\rm Hom}(A^n,A)\oplus \bigoplus_{k=1}^{n-1} {\rm
Hom}(C\otimes A^{n-k}, A\otimes C^k)\oplus {\rm Hom}(C,C^n).
$$
The cohomology of the complex $C_H((A,C)_\psi)$ is denoted
by $H_H((A,C)_\psi)$. 

The cohomology $H_H((A,C)_\psi)$ plays an important role in the
deformation theory of entwining structures. The latter can be developed
along the same lines as the deformation theory of algebra factorisations
in \cite{Brz:def}, following the Gerstenhaber deformation programme
\cite{Ger:def}. We sketch here the main results.

Let $(A,C)_\psi$ be an entwining
structure. A {\em formal deformation} of $(A,C)_\psi$ is an entwining structure
$(A_t,C_t)_{\psi_t}$, over the ring $k[[t]]$, where $A_t$, $C_t$ are
algebra  and coalgebra deformations of $A$
and $C$ respectively, and $\psi_t = \psi +\sum_{i=1}t^i\psi\sp{\su i}$,
$\psi\sp{\su i}:C\otimes A\to A\otimes C$. In other words the
deformation $(A_t,C_t)_{\psi_t}$ is characterised by three maps
expandable in the power series in $t$,
\begin{equation}
\mu_t = \mu +\sum_{i=1}t^i\mu\sp{\su i},\quad 
\Delta_t =\Delta +\sum_{i=1}t^i\Delta\sp{\su i}, \quad
\psi_t = \psi +\sum_{i=1}t^i\psi\sp{\su i}.
\label{expan}
\end{equation}
Two deformations $(A_t,C_t)_{\psi_t}$ and $(A_t,C_t)_{\tilde{\psi}_t}$
are {\em equivalent} to each other
 if there exist algebra isomorphism
$\alpha_t :A_t\to \tilde{A}_t$ and a coalgebra isomorphism 
$\gamma_t:C_t\to \tilde{C}_t$ 
of the form $\alpha_t = A +
\sum_{i=1} t^i\alpha^\su i$, $\gamma_t = C +
\sum_{i=1} t^i\gamma^\su i$, and such that $\tilde{\psi}_t\circ
(\gamma_t\otimes
\alpha_t) = (\alpha_t\otimes\gamma_t)\circ \psi_t$.
 A deformation
$(A_t,C_t)_{\psi_t}$ is called a {\em trivial deformation} if it is 
equivalent to an entwining structure over $k[[t]]$  
 in which all the maps $\mu^\su i$, $\Delta^\su
i$, $\psi^\su i$ in (\ref{expan}) vanish. 
An {\em infinitesimal deformation} of $(A,C)_\psi$ is a
deformation of $(A,C)_\psi$ modulo $t^2$. 

Using similar method to \cite[Theorem~3.1]{Brz:def} one proves the
following 
\begin{theorem}
There is a one-to-one correspondence between the equivalence classes of
infinitesimal deformations of $(A,C)_\psi$ 
and  $H^2_H((A,C)_\psi)$.
\label{deform.theorem}
\end{theorem}

Theorem~\ref{deform.theorem} is a standard result in various deformation
programmes. The appearance of the total cohomology of a double complex
makes the deformation theory of entwining structures similar to the
deformation theory of bialgebras \cite{GerSch:bia}. 

Furthermore, one can look  at the obstructions for extending a
deformation of an entwining structure modulo $t^n$ to a deformation
modulo $t^{n+1}$. Not surprisingly one finds that such obstructions are
classified by the third cohomology group $H^3_H((A,C)_\psi)$. We leave
the details of the analysis of obstructions to the reader, and only
mention that the details of a semi-dual case can be found in
\cite{Brz:def}.

\end{document}